\documentclass[12pt]{article}
\usepackage[reqno]{amsmath} 
\usepackage{amssymb, amsthm}
\usepackage{enumerate}  
\usepackage{hhline}
\usepackage{url} 
\usepackage{graphicx}
\usepackage{float}
\usepackage{nicefrac}
\usepackage{tikz}

\pgfdeclarelayer{edgelayer}
\pgfdeclarelayer{nodelayer}
\pgfsetlayers{edgelayer,nodelayer,main}

\usetikzlibrary{decorations.markings}

\tikzstyle{new}=[circle, minimum width=4pt,inner sep=0pt, fill=black,draw=black]
\tikzstyle{none}=[circle,fill=white,draw=black]
\tikzstyle{n}=[shape=rectangle,minimum width=1pt,inner sep=0pt, fill=none,draw=none]
\tikzstyle{emph}=[circle, minimum width=4pt,inner sep=0pt, fill=magenta,draw=magenta]

\tikzset{directed/.style={decoration={
  markings,
  mark=at position .6 with {\arrow{>}}},postaction={decorate}}}
\tikzset{nodirection/.style={}}

\newtheoremstyle{plainsl}%
	{\topsep}
	{\topsep}
	{\slshape} 
	{}
	{\normalfont\bfseries}
	{.}
	{ }
	{}

\swapnumbers

\theoremstyle{plainsl}
\newtheorem{theorem}{Theorem}[section]
\newtheorem{lemma}[theorem]{Lemma}
\newtheorem{conjecture}[theorem]{Conjecture}
\newtheorem{corollary}[theorem]{Corollary}
\newtheorem*{theorem*}{Theorem}

\renewcommand\proof{\noindent\textsl{Proof. }}
\newcommand\sqr[2]{{\vbox{\hrule height.#2pt
    \hbox{\vrule width.#2pt height#1pt \kern#1pt
        \vrule width.#2pt}\hrule height.#2pt}}}

\renewcommand\qed{%
	\ifmmode\eqno\sqr53
	\else\nolinebreak\ \hfill\sqr53\medbreak\fi}

\numberwithin{equation}{section}

\newcommand\mxm{\widehat{M}}
\newcommand\amx[1]{\widehat{M}(#1)}
\newcommand\amxb[2]{\widehat{M}_{#1}(#2)}
\newcommand\amxA[1]{\widehat{M}_A(#1)}
\newcommand\amxL[1]{\widehat{M}_L(#1)}
\newcommand\psd{\succcurlyeq}
\newcommand\cprod{\,\Box\,}

\newcommand\De{\Delta}

\newcommand\opk[1]{\mathop{\mathrm{#1}}\nolimits}
\newcommand\comp[1]{{\mkern2mu\overline{\mkern-2mu#1}}}

\newcommand\seq[3]{#1_{#2},\ldots,#1_{#3}}

\newcommand\pmat[1]{\begin{pmatrix} #1 \end{pmatrix}}

\newcommand\tr{\opk{tr}}

\title{Diagonal entries of the average mixing matrix}
\author{Chris Godsil\footnote{University of Waterloo, Waterloo, Canada. email: \protect\url{cgodsil@uwaterloo.ca}. C. Godsil gratefully acknowledges the support of the Natural Sciences and Engineering Council of Canada (NSERC), Grant No. RGPIN-9439.}, 
Krystal Guo\footnote{Universit\'{e} libre de Bruxelles, Brussels, Belgium. K. Guo gratefully acknowledges the support from ERC grant FOREFRONT
(grant agreement no. 615640) funded by the European Research Council under the EU’s
7th Framework Programme (FP7/2007-2013). Part of this research was done when K. Guo was a post-doctoral fellow at University of Waterloo. email: \protect\url{guo.krystal@gmail.com}}, 
Mariia Sobchuk\footnote{University of Waterloo, Waterloo, Canada. email: \protect\url{msobchuk@uwaterloo.ca}} }

\begin{document}

\maketitle

\begin{abstract}
	We study the diagonal entries of the average mixing matrix of continuous quantum walks. The average mixing matrix is a graph invariant; it is the sum of the Schur squares of spectral idempotents of the Hamiltonian. It is non-negative, doubly stochastic and positive semi-definite. We investigate the diagonal entries of this matrix.  We study the graphs for which the trace of the average mixing matrix is maximum or minimum and we classify those which are maximum. We give two constructions of graphs whose average mixing matrices have constant diagonal.
\end{abstract}

\section{Introduction}
\label{sec:intro}

Let $X$ be a graph on $n$ vertices with adjacency matrix $A$. We use $\De$ to denote
the $n\times n$ diagonal matrix with $\De_{i,i}$ equal to the valency of the vertex $i$.
The matrix $\De-A$ is the \textsl{Laplacian}, denoted $L(X)$, and $\De+A$ is known as 
the \textsl{signless Laplacian} of $X$.

A continuous random walk on a graph $X$ is given by the 1-parameter family of matrices
\[
	\exp(t(A-\De)).
\]
Physicists are concerned with quantum analogs. If $B$ is a linear combination of $A$ and $\De$,
the 1-parameter family of matrices
\[
	U(t) = \exp{itB}
\]
determines a \textsl{continuous quantum walk} on the vertices of $X$. We say that $U(t)$ is 
the \textsl{transition matrix} of the walk, and refer to $B$ as the \textsl{Hamiltonian}. 
(Usually $B$ is the adjacency matrix, the next most common 
choice is the Laplacian. If $X$ is regular, then the adjacency matrix and the two Laplacians provide 
the same information.) Kay \cite{KayReviewPST} provides a useful survey of continuous quantum walks.


Our overall aim is to establish connections between properties of a continuous quantum walk and properties
of the underlying graph. In this paper we focus on a matrix derived from $U(t)$, the
\textsl{mixing matrix} of the walk, which we define by
\[
	M(t) = U(t) \circ \comp{U(t)}.
\]
(Here $\circ$ denotes the Schur product of two matrices.) The \textsl{average mixing matrix}, denoted $\mxm$, 
is defined as follows:
\begin{equation}
	\label{eq:mxmdef}
	\amx{X} = \lim_{T\to\infty}\frac1T \int_0^T M(t)\,dt.
\end{equation}
The average mixing matrix is, intuitively, a distribution that the quantum walk adheres to, on average, over time, and thus may be thought of as a replacement for a stationary distribution. Earlier work on this matrix
appears in \cite{GodsilAverageMixing,TamonAdamzackAverageMixing, CoGdGuZh17, GoGuSi18}.

The average mixing matrix is doubly stochastic and positive semi-definite. In this paper, we investigate its diagonal entries. In particular, we study the graphs which attain the maximum and minimum of the trace of $\mxm$ with respect to the Laplacian and adjacency matrices.  We find that, amongst all graphs on $n$ vertices, $K_n$ maximizes the trace of average mixing matrix with respect to the Laplacian matrix. We also investigate the graphs whose average mixing matrices, with respect to the adjacency matrix, have constant diagonal, and we give a construction for non-regular graphs with such a property. Tables summarizinbg the results of some of our computations are provided.

\section{The Average Mixing Matrix}\label{sec:prelim}

Let $X$ be a graph on $n$ vertices and let $B \in \{A(X), L(X)\}$. Let $\theta_1, \ldots, \theta_d$ be the distinct eigenvalues of $B$ and, for $r=1,\ldots,d$, let $E_r$ be the idempotent projection onto the $\theta_r$ eigenspace of $B$; the spectral decomposition of $B$ is as follows:
\[
	B = \sum_{r=1}^d \theta_r E_r.
\]
The following is an important theorem, as it allows us to understand the average mixing matrix.

\begin{theorem}\cite{GodsilAverageMixing} \label{thm:godsilavgmix}
	Let $X$ be a graph and let $B \in \{A(X), L(X)\}$. Let $B = \sum_{r=0}^d \theta_r E_r$ 
	be the spectral decomposition of $B$. The average mixing matrix of $X$ with respect to $B$ is 
	\[
		\amx{B} = \sum_{r=0}^d E_r \circ E_r.\qed 
	\]
\end{theorem}

By way of example, we consider the complete graph $K_n$, with the adjacency matrix as Hamiltonian. Here $n-1$ is a simple eigenvalue and the asssociated eigenspace is spanned by the constant vectors; the projection is
\[
	E_1 = \frac1n J.
\]
The second eigenvalue of $K_n$ is $-1$, with multiplicity $n-1$. Since the spectral projections sum to $I$,
the second projection is
\[
	E_2 = I - \frac{1}{n}J.
\]
Hence the average mixing matrix of $K_n$ is
\[
	\frac1{n^2}J\circ J + \left(I - \frac{1}{n}J\right)\circ \left(I - \frac{1}{n}J\right)
		= \left(1-\frac2n\right)I +\frac2{n^2}J
\]

We will use $\amxb{B}{X}$ to denote the average mixing matrix of $X$, with the matrix $B$ as the Hamiltonian
of the quantum walk. We will denote by $\comp{X}$ the complement of a graph $X$. We will write $J_{j,k}$ 
for the $j\times k$ matrix with all entries equal to $1$; when $j = k,$ we will write $J_j$ for convenience. 
For Laplacian of connected graphs, the average mixing matrix is the same as complement, except when complement 
is not connected; we will give the explicit relation between $\amxb{L}{\comp{X}}$ and $\amxb{L}{X}$ 
for completeness, in Lemma \ref{lem:lapamx}.

\section{Laplacians}

Lemma \ref{prop:lap} summarizes some basic properties of Laplacian eigenvalues, which can be found in any 
standard algebraic graph theory text, including \cite{GR}.

\begin{lemma}\label{prop:lap} 
	Let $X$ be a graph on $n$ vertices and $L :=L(X)$ be the Laplacian matrix of $X$. Suppose the 
	distinct eigenvalues of $L$ are $\theta_0 \leq \cdots \leq \theta_d$, with corresponding 
	projections $E_0, \ldots, E_d$ onto the eigenspaces of $L$.
	\begin{enumerate}[(i)]
	\item For every $i$, we have that $n- \theta_i$ is an eigenvalues of $L(\comp{X})$.
	\item $\theta_0 = 0$ and the multiplicity of $0$ as an eigenvalue of $L$ is equal to the number 
	of connected components of $X$.
	\item The multiplicity of $n$ as an eigenvalue of $X$ is $c-1$, where $c$ is the number of 
	connected components of the complement of $X$.
	\item If $X$ has $c$ components $C_1, C_2, \ldots, C_c$ and the vertices of $X$ are 
	ordered $( x_{i_1}, \ldots, x_{i_{|C_i|}})_{i = 1}^{c}$, then
	\[
	E_0 = \pmat{\frac{1}{|C_1|}J_{|C_1|} & \mathbf{0} & \cdots & \mathbf{0} \\
				\mathbf{0} & \frac{1}{|C_2|}J_{|C_2|} &  	   & \mathbf{0} \\
				\mathbf{0} & \mathbf{0} & \ddots & \vdots \\
				\mathbf{0} & \mathbf{0} & \cdots & \frac{1}{|C_c|}J_{|C_c|} }
	\]
	where $\mathbf{0}$ denotes an all zero matrix of the appropriate order.
 \end{enumerate}
\end{lemma}

\begin{lemma}\label{lem:lapamx} Let $X$ be a connected  graph on $n$ vertices and $L :=L(X)$ be the Laplacian matrix of $X$. Suppose the distinct eigenvalues of $L$ are $\theta_0 \leq \cdots \leq \theta_d$, with corresponding projections $E_0, \ldots, E_d$ onto the eigenspaces of $L$.  If $\comp{X}$ has $c$ components $C_1, C_2, \ldots, C_c$ and the vertices of $X$ are ordered $( x_{i_1}, \ldots, x_{i_{|C_i|}})_{i = 1}^{c}$, then
\begin{equation}\label{eq:lapamx}
\amxL{\comp{X}} \! = \!\amxL{X} \!-\! \frac{1}{n^2} J_n \!+\!\pmat{\frac{1}{|C_1|^2}J_{|C_1|} & \mathbf{0} & \cdots & \mathbf{0} \\
			\mathbf{0} & \frac{1}{|C_2|^2}J_{|C_2|} &  	   & \mathbf{0} \\
			\mathbf{0} & \mathbf{0} & \ddots & \vdots \\
			\mathbf{0} & \mathbf{0} & \cdots & \frac{1}{|C_c|^2}J_{|C_c|} }
\end{equation}
where  $\mathbf{0}$ denotes an all zero matrix of the appropriate order.  \end{lemma}

\proof 
Observe that
\[
L(X) + L(\comp{X}) = -J + nI
\]
As in the statement of the lemma, let the distinct eigenvalues of $L$ be denoted by $\theta_0 \leq \cdots \leq \theta_d$, with corresponding projections $E_0, \ldots, E_d$ onto the eigenspaces of $L$.
If $i \neq 0$, then $L(\comp{X})E_i = (n- \theta_i)E_i$ and thus, since $\sum_{i =0 }^d E_i = I_n$
\begin{equation}\label{eq:lapeq2}
L(\comp{X}) = 0\cdot E_0 + \sum_{i =1}^d (n-\theta_i) E_i .
\end{equation}
If $\comp{X}$ is connected, then $\theta_i \neq n$ for any $i$ and thus \eqref{eq:lapeq2} is the spectral decomposition of $L(\comp{X})$.

For the more general result, we suppose that $\comp{X}$ has $c$ components $C_1, C_2, \ldots, C_c$ and the vertices of $X$ are ordered $( x_{i_1}, \ldots, x_{i_{|C_i|}})_{i = 1}^{c}$. Since $\theta_d = n$, we obtain that
\begin{equation}\label{eq:lapeq3}
	L(\comp{X}) = 0\cdot (E_0 + E_n) + \sum_{i =1}^{d-1} (n-\theta_i) E_i .
\end{equation}
We see that $E_0 + E_n$ is idempotent and thus \eqref{eq:lapeq3} is the spectral decomposition of $L(\comp{X})$ and the lemma follows.\qed

In particular, the Lemma \ref{lem:lapamx} implies that if $\comp{X}$ is connected, then $\amxL{X} = \amxL{\comp{X}}$.

\begin{corollary} 
	Let $X$ be a graph on $n$ vertices and $L :=L(X)$ be the Laplacian matrix of $X$. Suppose the 
	distinct eigenvalues of $L$ are $\theta_0 \leq \cdots \leq \theta_d$, with corresponding 
	projections $E_0, \ldots, E_d$ onto the eigenspaces of $L$.  If $\comp{X}$ has $c$ 
	components $C_1, C_2, \ldots, C_c$ , then
	\[
		\tr(\amxL{\comp{X}}) = \tr(\amxL{X}) - \frac{1}{n} + \sum_{i = 1}^c \frac{1}{|C_i|}.
	\]
\end{corollary}

\proof
This follows from taking the trace of both sides of \eqref{eq:lapamx}.\qed

\section{An Ordering}

For symmetric, square matrices $A$ and $B$, we write $A\psd B$ if $A-B$ is positive semidefinite. 
The relation $\psd$ is a useful partial ordering on matrices. We briefly investigate
some of its properties when applied to average mixing matrices.

\begin{lemma}\label{lem:ord}
	Let $X$ and $Y$ be graphs on the same vertex set. If each spectral idempotent of
	$A(Y)$ is the sum of spectral idempotents of $A(X)$, then $\amxA{Y}\psd\amxA{X}$.  
	Similarly, if each spectral idempotent of
	$L(Y)$ is the sum of spectral idempotents of $L(X)$, then $\amxL{Y}\psd\amxL{X}$.
\end{lemma}

\proof
We have
\[
	(E+F)^{\circ2} = E^{\circ2} + F^{\circ2} + 2E\circ F.
\]
If $E$ and $F$ are positive semidefinite, so are the three terms in the sum above,
whence
\[
	(E+F)^{\circ2} \psd E^{\circ2} + F^{\circ2}.
\]
We apply this iteratively and, using Theorem \ref{thm:godsilavgmix}, the lemma follows.\qed

Lemma \ref{lem:ord} has many consequences, some of which will be explored in the next section. For now, give a lemma about the products of graphs, which uses the same basic idea in the proof. Recall that  $X\cprod Y$ denotes the Cartesian product of graphs $X$ and $Y$ and $X\times Y$ denotes the categorical (or direct) product of $X$ and $Y$.

\begin{lemma}
	Let  $X$ and $Y$ be graphs.
	\begin{enumerate}[(a)]
	\item $\{\amxA{X\cprod Y}, \amxA{X\times Y}\} \psd \amxA{X}\otimes\amxA{Y}$; and
	\item $\{\amxL{X\cprod Y}, \amxA{X\cprod Y}\} \psd \amxL{X}\otimes\amxL{Y}$.
	\end{enumerate}
\end{lemma}

\proof 
Let $\seq E1\ell$ and $\seq F1k$ be the respective spectral idempotents of $X$ and $Y$.
Then each spectral idempotent of $X\cprod Y$ is a sum of idempotents of the form
$E_r\otimes F_s$. To complete the proof that $\amxA{X\cprod Y} \psd \amxA{X}\otimes\amxA{Y}$, note that
\[
	(E_r\otimes F_s)^{\circ2} = E_r^{\circ2} \otimes F_s^{\circ2}.
\]
This argument will also hold when we use the Laplacian in place of the adjacency matrix.
The proof of part (b) also follows similarly.\qed

If $P$ is a permutation matrix and $A(Y)=P^TA(X)P$, then
\[
	\amx{Y} = P^T\amx{X} P.
\]
This indicates that the partial ordering on average mixing matrices using $\psd$ cannot
generally correspond to a useful ordering on graphs.

\section{Application to Quantum Walks using the Laplacian}\label{sec:lapqw}

We will use Lemma \ref{lem:ord} to show that $K_n$ maximizes the trace of average mixing matrix with 
respect to the Laplacian matrix, amongst all graphs on $n$ vertices.

Note that, for graph $X$ and $Y$, we denote by $X+Y$ the disjoint union of $X$ and $Y$.

\begin{lemma}\label{lem:ordlap} 
	Let $X$ be a connected graph on $n$ vertices. Then \[ \amxL{K_n} \psd \amxL{X}.
	\]
	More generally, if $X$ has $m$ connected components $C_1,\ldots,C_m$ with $c_1,\ldots,c_m$ 
	vertices respectively, then
	\[
		\amxL{K_{c_1} + \cdots + K_{c_m}} \psd \amxL{X}.
	\]
\end{lemma}

\proof 
Let $I$ and $J$ be the $n\times n$ identity and all ones matrix, respectively. We observe that
\begin{equation} \label{eq:lapkn} 
	L(K_n) = nI - J = 0\left( \frac{1}{n} J \right) + n \left( I - \frac{1}{n} J \right).  
\end{equation}
We see that the right side of \eqref{eq:lapkn} gives the spectral decomposition of $L(K_n)$.

Any connected graph $X$ on $n$ vertices has $\frac{1}{n}J$ as a spectral idempotent and so the sum of the other spectral idempotents of $X$ must be equal to $ I - \frac{1}{n} J$. Thus  we have obtained that each spectral idempotent of $L(K_n)$ is the sum of spectral idempotents of $L(X)$, and so $\amxL{Y}\psd\amxL{X}$ by Lemma \ref{lem:ord}. The more general statement follows from consider the connected components separately and recalling that, for graphs $X $ and $Y$,
\[
	\amxL{X+Y} = \pmat{\amxL{X} & 0 \\0 & \amxL{Y}}
\]
and the lemma follows.\qed

\begin{corollary}\label{lem:laptrs}
	Let $X$ be a connected graph on $n$ vertices. Then 
	\[
		\tr(\amxL{K_n}) \geq \tr(\amxL{X}),
	\]
	and equality holds if and only if $X$ is isomorphic to $K_n$.
	More generally, if $X$ has $m$ connected components $C_1,\ldots,C_m$ with $c_1,\ldots,c_m$ 
	vertices respectively, then
	\[\tr(\amxL{K_{c_1} + \cdots + K_{c_m}}) \geq \tr(\amxL{X}).\]
	Equality hold if and only if $X$ is isomorphic to $K_{c_1} + \cdots + K_{c_m}$.
\end{corollary}

\proof 
Recall the following property from elementary linear algebra: if $A \psd B$, then $\tr(A) \geq \tr(B)$. We immediately obtain for any connected graph $X$ on $n$ vertices that
\[
\tr(\amxL{K_n}) \geq \tr(\amxL{X}).
\]
Note that for positive semidefinite matrices $A$ and $B$, if $A \psd B$ and $\tr(A) = \tr(B)$, 
then $\tr(A-B) = 0$ and $A-B\psd 0$ and so $A-B = 0$.  Thus, if $\tr(\amxL{K_n}) =\tr(\amxL{X})$ 
then $\amxL{K_n} =\amxL{X}$ and thus $X$ is isomorphic to $K_n$. [maybe some more argument is needed here]

The more general statement follows if we consider the connected components separately.\qed

From \eqref{eq:lapkn}, we can see that
\[ \amxL{K_n} =  \left( \frac{1}{n} J \right)^{\circ 2} + \left( I - \frac{1}{n} J \right) ^{\circ 2}= \left( 1 - \frac{2}{n}  \right)I +\frac{2}{n^2} J
\]
and so
\[
\tr(\amxL{K_n}) = n - 2 + \frac{2}{n}.
\]
Note that $\tr(\amxL{\comp{K_n}}) = \tr(I) = n$.

\section{Eigenspace Refinement}\label{sec:adjqw}

Similar to Section \ref{sec:lapqw}, we look at implications of Lemma \ref{lem:ord} for the average mixing 
matrix with respect to the adjacency matrix.

\begin{lemma} 
	If $X$ has an equitable partition with parts of size $a_1,\ldots, a_m$, then the eigenspaces 
	of $X$ refine the eigenspaces of the disjoint union of $K_{a_1}, \ldots, K_{a_m}$.  
\end{lemma}

\proof 
This follows because every eigenvector of $X$ is either constant, or sums to $0$ on each part of the equitable partition.\qed

The trace of $\amx{A(K_n)}$ is
\[
\frac{1+(n-1)^2}{n} = \frac{n^2-2n + 2}{n}.
\]
Thus, with Lemma \ref{lem:ord}, we obtain the following:

\begin{lemma} 
	If $X$ has an equitable partition with parts of size $a_1,\ldots, a_m$, then
	\[ \
		tr{\amxA{X}} \leq \sum_{j=1}^m  \frac{a_j^2-2a_j + 2}{a_j}
	\]  
\end{lemma}

\proof 
Let $Y = K_{a_1} \cup \ldots \cup K_{a_m}$, the disjoint union of complete graphs. By the previous lemma, 
we see that the eigenspaces of $X$ refine those of $Y$ and thus
\[
	\tr{\amxA{X}} \leq \tr{\amxA{Y}}.
\]
We see that
\[
	\tr(\widehat{M}(Y)) = \sum_{r=1}^m \tr(\widehat{M}(K_{a_r}))=  \sum_{r=1}^m  \frac{a_r^2-2a_r + 2}{a_r}
\]
and the result follows.\qed

We see that $a_1 + \cdots + a_m = n$. Thus
\[ 
\begin{split}
	\sum_{j=1}^m  \frac{a_j^2 - 2a_j + 2}{a_j} &= \sum_{j=1}^m \left( a_j -2 +  \frac{2}{a_j}\right) \\
	&= n - 2m  + 2 \sum_{j=1}^m \frac{1}{a_j} \\
\end{split}
\]

\begin{corollary}\label{cor:regKnMaxTr} 
	If $X$ is a regular graph on $n$ vertices, then \[\tr{\amxA{X}} \leq \tr{\amxA{K_n}}\].
\end{corollary}

\section{Constant Diagonal}\label{sec:constdiag}

For this section, we will use $\amx{X}$ to denote $\amxA{X}$.

We consider graphs whose average mixing matrices with respect to the adjacency matrix have 
constant diagonal. These matrices are easy to work with because of the following trivial lemma. 

\begin{lemma} 
	Let $M$ be a positive semi-definite matrix. If $M$ has a constant diagonal, 
	then $M_{u,u} \geq M_{u,v}$. Furthermore, if $M$ has a constant diagonal 
	and $M_{u,u} = M_{u,v}$, then $Me_u = Me_v$. 
\end{lemma}

\proof
A $2\times 2$ principal submatrix of $M$ is the Gram matrix of two vectors, in our case of the same length.
Hence the result is a slightly disguised version of Cauchy-Schwarz.\qed

A graph is \textsl{walk-regular} if the number of closed walks at vertex $v$ of length $k$ is $c_k$, 
a constant independent of the choice of vertex, for every $k$. Equivalently, $X$ is walk-regular $A(X)^k$ 
has a constant diagonal for every $k$. We refer to \cite{GR} for further details.

\begin{lemma} 
	If $X$ is a walk-regular graph, then $\amx{X}$ has a constant diagonal. 
\end{lemma}

\proof 
Since the eigenprojections of $A(X)$ are polynomials in $A(X)$, we get that they are matrices with constant diagonals. Thus, $\amx{S}$ is a sum of matrices with constant diagonals.\qed

It is surprising that these are not the only graphs who average mixing matrices have constant diagonals. 
We will consider rooted products of graphs.

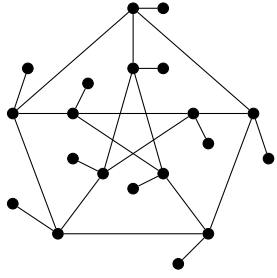
\begin{figure}
\centering
\begin{tikzpicture}[scale=0.8]
	\begin{pgfonlayer}{nodelayer}
		\node [style=new] (0) at (-1, 0) {};
		\node [style=new] (1) at (1, 0) {};
		\node [style=new] (2) at (-0.5, -1) {};
		\node [style=new] (3) at (0.5, -1) {};
		\node [style=new] (4) at (0, 0.75) {};
		\node [style=new] (5) at (0, 1.75) {};
		\node [style=new] (6) at (-2, 0) {};
		\node [style=new] (7) at (-1.25, -2) {};
		\node [style=new] (8) at (1.25, -2) {};
		\node [style=new] (9) at (2, 0) {};
		\node [style=new] (10) at (0.5, 1.75) {};
		\node [style=new] (11) at (2.25, -0.75) {};
		\node [style=new] (12) at (0.75, -2.5) {};
		\node [style=new] (13) at (-2, -1.5) {};
		\node [style=new] (14) at (-1.75, 0.75) {};
		\node [style=new] (15) at (0.5, 0.75) {};
		\node [style=new] (16) at (1.25, -0.5) {};
		\node [style=new] (17) at (0, -1.25) {};
		\node [style=new] (18) at (-1, -0.75) {};
		\node [style=new] (19) at (-0.75, 0.5) {};
	\end{pgfonlayer}
	\begin{pgfonlayer}{edgelayer}
		\draw [style=nodirection] (0) to (1);
		\draw [style=nodirection] (1) to (2);
		\draw [style=nodirection] (2) to (4);
		\draw [style=nodirection] (4) to (3);
		\draw [style=nodirection] (3) to (0);
		\draw [style=nodirection] (0) to (6);
		\draw [style=nodirection] (6) to (7);
		\draw [style=nodirection] (7) to (8);
		\draw [style=nodirection] (8) to (9);
		\draw [style=nodirection] (9) to (5);
		\draw [style=nodirection] (5) to (6);
		\draw [style=nodirection] (5) to (4);
		\draw [style=nodirection] (2) to (7);
		\draw [style=nodirection] (3) to (8);
		\draw [style=nodirection] (1) to (9);
		\draw [style=nodirection] (10) to (5);
		\draw [style=nodirection] (14) to (6);
		\draw [style=nodirection] (13) to (7);
		\draw [style=nodirection] (12) to (8);
		\draw [style=nodirection] (11) to (9);
		\draw [style=nodirection] (15) to (4);
		\draw [style=nodirection] (19) to (0);
		\draw [style=nodirection] (18) to (2);
		\draw [style=nodirection] (17) to (3);
		\draw [style=nodirection] (16) to (1);
	\end{pgfonlayer}
\end{tikzpicture}
\caption{The rooted product of the Petersen graph with $K_2$.\label{fig:rootedPete} }
\end{figure}

Let $X$ be a graph with vertices $\{v_1,\ldots, v_n\}$ and let $Y$ be a disjoint union of rooted graphs $Y_1, \ldots, Y_n$, rooted at $y_1,\ldots, y_n$, respectively. The \textsl{rooted product} of $X$ and $Y$, denoted $X(Y)$, is the graph obtain by identifying $v_i$ with the root vertex of $Y_i$. The rooted product was first introduced by Godsil and McKay in \cite{GoMc78}. We will consider the special case where $Y$ is a sequence of $n$ copies of $K_2$. In this case, we will write $X(K_2)$ to denote the rooted product. Figure \ref{fig:rootedPete} showed the rooted product of the Petersen graph with $K_2$, which is not walk-regular but whose average mixing matrix has a constant diagonal, by Corollary \ref{cor:rootedmxm}.

In \cite{GoGuSi18}, the authors consider the rooted product of $X$ with $K_2$. The results Lemma~3.3
and Theorem~3.4 there are stated under the assumption that ther eigenvalues of $X$ are simple,
but the proofs make no use of this assumption. We state the results for general graphs here and provide a brief proof.

\begin{lemma}\label{lem:rootedidems}
	Let $X$ be a graph and let $F_1,\ldots,F_d$ the orthogonal projections onto the eigenspaces of $A(X)$ 
	with corresponding eigenvalues $\lambda_1,\ldots, \lambda_d$.
	Then, the eigenvalues of the adjacency matrix of $X(K_2)$ are $\{\mu_i,\nu_i\}_{i=1,d}$ 
	where $\mu_i, \nu_i$ are the roots of $t^2 - \lambda_i t - 1 = 0$. For $\mu \in \{\mu_i, \nu_i\}$, 
	the projection onto the $\mu$-eigenspace is
	\[
		\frac{1}{\mu^2 + 1} \pmat{\mu^2 F_i & \mu F_i \\ \mu F_i & F_i}.
	\]
\end{lemma}

\proof 
The eigenvalues of $X(K_2)$ are given in \cite{GoMc78}. Fix $i$ and $\mu \in \{\mu_i, \nu_i\}$.
\[
E := \frac{1}{\mu^2 + 1} \pmat{\mu^2 F_i & \mu F_i \\ \mu F_i & F_i}.
\]
Observe the $E^2 = E$ and
\[
\begin{split}
A\left(X(K_2)\right)E  &=  \frac{1}{\mu^2 + 1} \pmat{ A & I \\ I & 0}\pmat{\mu^2 F_i & \mu F_i \\ \mu F_i & F_i} \\
&= \frac{1}{\mu^2 + 1}\pmat{(\lambda\mu^2 + \mu) F_i & (\lambda\mu + 1) F_i \\ \mu^2 F_i & \mu F_i} \\
&= \frac{\mu}{\mu^2 + 1}\pmat{(\lambda\mu + 1) F_i & \frac{(\lambda\mu + 1)}{\mu} F_i \\ \mu F_i & F_i}.
\end{split}
\]
Recall that $\mu^2 - \lambda \mu - 1 =0$ and so $\lambda \mu + 1 = \mu^2$. We obtain that $A\left(X(K_2)\right)E  = \mu E$ and the lemma follows.\qed

\begin{table}[htb]
\begin{tabular}{l|l|l|l}
$n$ & number of graphs & number of graphs  & number of walk regular \\
& on $n$ where $\amxA{X}$ has & on $n$ where $\amxL{X}$ has  & graphs on $n$ vertices\\
& constant diagonal & constant diagonal &  \\
\hline
2 & 2     & 2                  & 2              \\
3 & 2      & 2                 & 2              \\
4 & 7       & 5               & 4              \\
5 & 3        & 3               & 3              \\
6 & 15        & 12               & 8              \\
7 & 4         & 4             & 4              \\
8 & 48        & 59             & 14             \\
9 & 12       & 9              & 9
\end{tabular}
\caption{Graphs whose average mixing matrix has a constant diagonal with respect to the  adjacency matrix and Laplacian adjacency matrix. \label{tab:Aconstdiag}}
\end{table}

\begin{theorem}\cite{GoGuSi18}\label{thm:rootedmhat}
	Let $X$ be a graph and let $F_1,\ldots,F_d$ the orthogonal projections onto the eigenspaces of $A(X)$ 
	with corresponding eigenvalues $\lambda_1,\ldots, \lambda_d$.
	Then
	\begin{equation*}
		\mxm(X(K_2)) = \pmat{\mxm(X) - N & N \\ N & \mxm(X) - N}
	\end{equation*}
	where
	\begin{equation*} 
		N =  \sum_{i=1}^n  \left(\frac{2}{\lambda_i^2 + 4}\right)(F_i\circ F_i).
	\end{equation*}
\end{theorem}

\begin{corollary}\label{cor:rootedmxm} 
	If $X$ is a walk-regular graph and $Y$ is the rooted graph of $X$ with $K_2$, then $\amx{X}$ has 
	constant diagonal. 
\end{corollary}

\proof 
Note that for a walk-regular graph $X$, matrices $\mxm(X)$ and $N$ both have constant diagonal, and thus 
the corollary follows.\qed

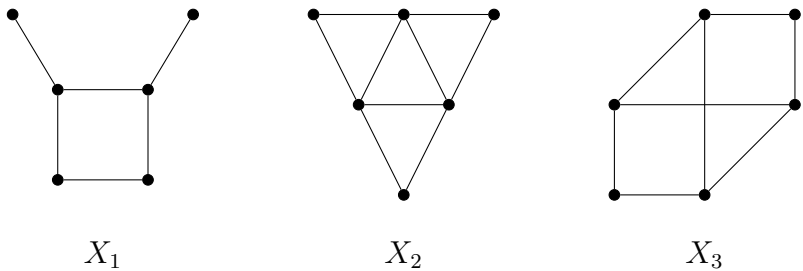
\begin{figure}[tbp]
\centering
\begin{tikzpicture}[scale=0.8]
	\begin{pgfonlayer}{nodelayer}
		\node [style=new] (0) at (-1.5, 1) {};
		\node [style=new] (1) at (0, 1) {};
		\node [style=new] (2) at (1.5, 1) {};
		\node [style=new] (3) at (0, -2) {};
		\node [style=new] (4) at (0.75, -0.5) {};
		\node [style=new] (5) at (-0.75, -0.5) {};
		\node [style=new] (6) at (-5.75, -0.25) {};
		\node [style=new] (7) at (-4.25, -0.25) {};
		\node [style=new] (8) at (-5.75, -1.75) {};
		\node [style=new] (9) at (-4.25, -1.75) {};
		\node [style=new] (10) at (-6.5, 1) {};
		\node [style=new] (11) at (-3.5, 1) {};
		\node [style=new] (12) at (3.5, -2) {};
		\node [style=new] (13) at (6.5, 1) {};
		\node [style=new] (14) at (3.5, -0.5) {};
		\node [style=new] (15) at (5, 1) {};
		\node [style=new] (16) at (5, -2) {};
		\node [style=new] (17) at (6.5, -0.5) {};
		\node [style=n] (18) at (-5, -3) {$X_1$};
		\node [style=n] (19) at (-5, -2.75) {};
		\node [style=n] (20) at (0, -3) {$X_2$};
		\node [style=n] (21) at (5, -3) {$X_3$};
	\end{pgfonlayer}
	\begin{pgfonlayer}{edgelayer}
		\draw [style=nodirection] (0) to (5);
		\draw [style=nodirection] (5) to (3);
		\draw [style=nodirection] (3) to (4);
		\draw [style=nodirection] (4) to (5);
		\draw [style=nodirection] (5) to (1);
		\draw [style=nodirection] (1) to (4);
		\draw [style=nodirection] (4) to (2);
		\draw [style=nodirection] (2) to (1);
		\draw [style=nodirection] (1) to (0);
		\draw [style=nodirection] (10) to (6);
		\draw [style=nodirection] (6) to (8);
		\draw [style=nodirection] (8) to (9);
		\draw [style=nodirection] (9) to (7);
		\draw [style=nodirection] (7) to (6);
		\draw [style=nodirection] (7) to (11);
		\draw [style=nodirection] (12) to (14);
		\draw [style=nodirection] (12) to (16);
		\draw [style=nodirection] (16) to (17);
		\draw [style=nodirection] (17) to (14);
		\draw [style=nodirection] (14) to (15);
		\draw [style=nodirection] (15) to (16);
		\draw [style=nodirection] (15) to (13);
		\draw [style=nodirection] (13) to (17);
	\end{pgfonlayer}
\end{tikzpicture}
\caption{Graphs $X1, X_2, X_3$ on $6$ vertices such that $\mxm(X_i)$ has constant diagonal.\label{fig:constdiagNotReg} }
\end{figure}

Two sources of graphs whose average mixing matrices have constant diagonals are walk-regular graphs and their rooted product with $K_2$. A computation on graphs up to $9$ vertices reveals that there are other graphs with the property. Figure \ref{fig:constdiagNotReg} shows three graphs $X1, X_2, X_3$ on $6$ vertices such that $\mxm(X_i)$ has constant diagonal, but $X_i$ is not a walk-regular graph nor a walk-regular graph rooted with $K_2$. Table \ref{tab:Aconstdiag} shows the numbers of graphs whose average mixing matrix has a constant diagonal with respect to the  adjacency matrix and the Laplacian matrix.

\section{Trace computations}

\begin{table}[hb]
\centering
\begin{tabular}{l|c|c|c|c|c|c}
\hline
$n$ & $3$ & $4$ & $5$ & $6$ & $7$ & $8$ \\
\hline
$\max_X \tr(\amxA{X})$ & $1\nicefrac{2}{3}$ & $2\nicefrac{1}{2}$ & $3\nicefrac{2}{5}$ & $4\nicefrac{1}{3}$ & $5\nicefrac{2}{7}$ & $6\nicefrac{1}{4}$ \\
\hline
Graphs & $K_3$  & $K_4$ & $K_5$ & $K_6$ & $K_7$ & $K_8$\\
\hline
\end{tabular}
\caption{Graphs on $n$ vertices attaining the maximum trace with respect to $\amxA$ for $n = 3,4,5,6,7,8$.\label{tab:Amaxtr}}
\end{table}

Tables for $A$ and $L$ of maximum and minimum trace and the graphs which attain the maximum and minimum. We have determined the graphs attaining the maximum trace with respect to the Laplacian matrix. Tables \ref{tab:Amaxtr}, \ref{tab:Lmintr}, and \ref{tab:Amintr}  show the minimum with respect to $L$ and the min and max with respect to $A$.

\begin{table}[tbp]
\centering
\begin{tabular}{|l|l|c|}
\hline
$n$ & $\min_X \tr(\amxL{X})$ & Graphs \\
\hline
$3$ & $1\nicefrac{1}{3}$ & \begin{tikzpicture}[scale=0.6]
	\begin{pgfonlayer}{nodelayer}
		\node [style=new] (0) at (0, 0) {};
		\node [style=new] (1) at (-1, 0) {};
		\node [style=new] (2) at (1, 0) {};
	\end{pgfonlayer}
	\begin{pgfonlayer}{edgelayer}
		\draw [style=nodirection] (1) to (0);
		\draw [style=nodirection] (0) to (2);
	\end{pgfonlayer}
\end{tikzpicture} \\
\hline
$4$ & $1\nicefrac{1}{4}$ & \begin{tikzpicture}[scale=0.6]
	\begin{pgfonlayer}{nodelayer}
		\node [style=new] (0) at (0, 0) {};
		\node [style=new] (1) at (-1, 0) {};
		\node [style=new] (2) at (1, 0) {};
		\node [style=new] (3) at (-2, 0) {};
	\end{pgfonlayer}
	\begin{pgfonlayer}{edgelayer}
		\draw [style=nodirection] (1) to (0);
		\draw [style=nodirection] (0) to (2);
		\draw [style=nodirection] (3) to (1);
	\end{pgfonlayer}
\end{tikzpicture} \\
\hline
$5$ & $1\nicefrac{2}{5}$ & \begin{tikzpicture}[scale=0.6]
	\begin{pgfonlayer}{nodelayer}
		\node [style=new] (0) at (0, 0) {};
		\node [style=new] (1) at (-1, 0) {};
		\node [style=new] (2) at (1, 0) {};
		\node [style=new] (3) at (-2, 0) {};
		\node [style=new] (4) at (-3, 0) {};
		\node [style=new] (5) at (3, -1) {};
		\node [style=new] (6) at (3, 0) {};
		\node [style=new] (7) at (4, 1) {};
		\node [style=new] (8) at (5, 0) {};
		\node [style=new] (9) at (5, -1) {};
	\end{pgfonlayer}
	\begin{pgfonlayer}{edgelayer}
		\draw [style=nodirection] (1) to (0);
		\draw [style=nodirection] (0) to (2);
		\draw [style=nodirection] (3) to (1);
		\draw [style=nodirection] (4) to (3);
		\draw [style=nodirection] (6) to (5);
		\draw [style=nodirection] (5) to (9);
		\draw [style=nodirection] (9) to (8);
		\draw [style=nodirection] (8) to (7);
		\draw [style=nodirection] (7) to (6);
		\draw [style=nodirection] (6) to (8);
	\end{pgfonlayer}
\end{tikzpicture} \\
\hline
$6$ & $1\nicefrac{1}{3}$ & \begin{tikzpicture}[scale=0.6]
	\begin{pgfonlayer}{nodelayer}
		\node [style=new] (0) at (0, 0) {};
		\node [style=new] (1) at (-1, 0) {};
		\node [style=new] (2) at (1, 0) {};
		\node [style=new] (3) at (-2, 0) {};
		\node [style=new] (4) at (-3, 0) {};
		\node [style=new] (5) at (3, -1) {};
		\node [style=new] (6) at (3, 0) {};
		\node [style=new] (7) at (4, 1) {};
		\node [style=new] (8) at (5, 0) {};
		\node [style=new] (9) at (5, -1) {};
		\node [style=new] (10) at (-4, 0) {};
		\node [style=new] (11) at (4, -2) {};
	\end{pgfonlayer}
	\begin{pgfonlayer}{edgelayer}
		\draw [style=nodirection] (1) to (0);
		\draw [style=nodirection] (0) to (2);
		\draw [style=nodirection] (3) to (1);
		\draw [style=nodirection] (4) to (3);
		\draw [style=nodirection] (6) to (5);
		\draw [style=nodirection] (5) to (9);
		\draw [style=nodirection] (9) to (8);
		\draw [style=nodirection] (8) to (7);
		\draw [style=nodirection] (7) to (6);
		\draw [style=nodirection] (6) to (8);
		\draw [style=nodirection] (10) to (4);
		\draw [style=nodirection] (9) to (11);
		\draw [style=nodirection] (11) to (5);
		\draw [style=nodirection] (5) to (7);
		\draw [style=nodirection] (6) to (11);
	\end{pgfonlayer}
\end{tikzpicture}\\
\hline
$7$ & $1\nicefrac{3}{7}$ & \begin{tikzpicture}[scale=0.6]
	\begin{pgfonlayer}{nodelayer}
		\node [style=new] (0) at (0, 0) {};
		\node [style=new] (1) at (-1, 0) {};
		\node [style=new] (2) at (1, 0) {};
		\node [style=new] (3) at (-2, 0) {};
		\node [style=new] (4) at (-3, 0) {};
		\node [style=new] (5) at (-4, 0) {};
		\node [style=new] (6) at (-5, 0) {};
		\node [style=new] (7) at (3, 0) {};
		\node [style=new] (8) at (4, 0.75) {};
		\node [style=new] (9) at (4, -0.75) {};
		\node [style=new] (10) at (5, -0.75) {};
		\node [style=new] (11) at (6, -0.75) {};
		\node [style=new] (12) at (6, 0.75) {};
		\node [style=new] (13) at (5, 0.75) {};
		\node [style=new] (14) at (-4.5, -3) {};
		\node [style=new] (15) at (-3.5, -2.25) {};
		\node [style=new] (16) at (-3.5, -3.75) {};
		\node [style=new] (17) at (-2, -3.75) {};
		\node [style=new] (18) at (-2, -2.25) {};
		\node [style=new] (19) at (-1, -3) {};
		\node [style=new] (20) at (-2.75, -3) {};
		\node [style=new] (21) at (1.25, -2.25) {};
		\node [style=new] (22) at (1.25, -3.75) {};
		\node [style=new] (23) at (3.25, -4) {};
		\node [style=new] (24) at (3.25, -2) {};
		\node [style=new] (25) at (4.25, -3) {};
		\node [style=new] (26) at (2.25, -2.5) {};
		\node [style=new] (27) at (2.25, -3.5) {};
	\end{pgfonlayer}
	\begin{pgfonlayer}{edgelayer}
		\draw [style=nodirection] (1) to (0);
		\draw [style=nodirection] (0) to (2);
		\draw [style=nodirection] (3) to (1);
		\draw [style=nodirection] (4) to (3);
		\draw [style=nodirection] (5) to (4);
		\draw [style=nodirection] (6) to (5);
		\draw [style=nodirection] (7) to (8);
		\draw [style=nodirection] (8) to (13);
		\draw [style=nodirection] (13) to (10);
		\draw [style=nodirection] (10) to (11);
		\draw [style=nodirection] (11) to (12);
		\draw [style=nodirection] (12) to (13);
		\draw [style=nodirection] (10) to (9);
		\draw [style=nodirection] (9) to (8);
		\draw [style=nodirection] (9) to (7);
		\draw [style=nodirection] (14) to (15);
		\draw [style=nodirection] (15) to (18);
		\draw [style=nodirection] (18) to (19);
		\draw [style=nodirection] (19) to (17);
		\draw [style=nodirection] (17) to (16);
		\draw [style=nodirection] (16) to (14);
		\draw [style=nodirection] (15) to (16);
		\draw [style=nodirection] (17) to (18);
		\draw [style=nodirection] (14) to (20);
		\draw [style=nodirection] (20) to (19);
		\draw [style=nodirection] (16) to (20);
		\draw [style=nodirection] (20) to (18);
		\draw [style=nodirection] (21) to (26);
		\draw [style=nodirection] (26) to (27);
		\draw [style=nodirection] (27) to (21);
		\draw [style=nodirection] (21) to (22);
		\draw [style=nodirection] (22) to (26);
		\draw [style=nodirection] (27) to (22);
		\draw [style=nodirection] (21) to (24);
		\draw [style=nodirection] (24) to (25);
		\draw [style=nodirection] (25) to (23);
		\draw [style=nodirection] (23) to (22);
		\draw [style=nodirection] (27) to (23);
		\draw [style=nodirection] (23) to (24);
		\draw [style=nodirection] (24) to (26);
		\draw [style=nodirection] (26) to (25);
		\draw [style=nodirection] (25) to (27);
	\end{pgfonlayer}
\end{tikzpicture} \\
\hline
$8$ & $1\nicefrac{3}{8}$ & \begin{tikzpicture}[scale=0.6]
	\begin{pgfonlayer}{nodelayer}
		\node [style=new] (0) at (0, 0) {};
		\node [style=new] (1) at (-1, 0) {};
		\node [style=new] (2) at (1, 0) {};
		\node [style=new] (3) at (-2, 0) {};
		\node [style=new] (4) at (-3, 0) {};
		\node [style=new] (5) at (-4, 0) {};
		\node [style=new] (6) at (-5, 0) {};
		\node [style=new] (7) at (-6, 0) {};
		\node [style=new] (8) at (4.25, 1.25) {};
		\node [style=new] (9) at (5.25, 1.25) {};
		\node [style=new] (10) at (4.25, -1.25) {};
		\node [style=new] (11) at (5.25, -1.25) {};
		\node [style=new] (12) at (3.5, 0.5) {};
		\node [style=new] (13) at (3.5, -0.5) {};
		\node [style=new] (14) at (6, 0.5) {};
		\node [style=new] (15) at (6, -0.5) {};
		\node [style=new] (16) at (1.75, -2.75) {};
		\node [style=new] (17) at (1.75, -3.75) {};
		\node [style=new] (18) at (2.75, -2.75) {};
		\node [style=new] (19) at (2.75, -3.75) {};
		\node [style=new] (20) at (1, -2) {};
		\node [style=new] (21) at (3.5, -2) {};
		\node [style=new] (22) at (3.5, -4.5) {};
		\node [style=new] (23) at (1, -4.5) {};
		\node [style=new] (24) at (-2.25, -2) {};
		\node [style=new] (25) at (-2.25, -4.5) {};
		\node [style=new] (26) at (-3.25, -4.5) {};
		\node [style=new] (27) at (-4, -2.75) {};
		\node [style=new] (28) at (-4, -3.75) {};
		\node [style=new] (29) at (-3.25, -2) {};
		\node [style=new] (30) at (-1.5, -3.75) {};
		\node [style=new] (31) at (-1.5, -2.75) {};
	\end{pgfonlayer}
	\begin{pgfonlayer}{edgelayer}
		\draw [style=nodirection] (1) to (0);
		\draw [style=nodirection] (0) to (2);
		\draw [style=nodirection] (3) to (1);
		\draw [style=nodirection] (4) to (3);
		\draw [style=nodirection] (5) to (4);
		\draw [style=nodirection] (6) to (5);
		\draw [style=nodirection] (7) to (6);
		\draw [style=nodirection] (9) to (14);
		\draw [style=nodirection] (14) to (15);
		\draw [style=nodirection] (15) to (11);
		\draw [style=nodirection] (11) to (10);
		\draw [style=nodirection] (10) to (13);
		\draw [style=nodirection] (13) to (12);
		\draw [style=nodirection] (12) to (8);
		\draw [style=nodirection] (9) to (15);
		\draw [style=nodirection] (15) to (10);
		\draw [style=nodirection] (10) to (12);
		\draw [style=nodirection] (12) to (9);
		\draw [style=nodirection] (8) to (14);
		\draw [style=nodirection] (14) to (11);
		\draw [style=nodirection] (11) to (13);
		\draw [style=nodirection] (13) to (8);
		\draw [style=nodirection] (20) to (21);
		\draw [style=nodirection] (21) to (22);
		\draw [style=nodirection] (22) to (23);
		\draw [style=nodirection] (23) to (20);
		\draw [style=nodirection] (16) to (17);
		\draw [style=nodirection] (17) to (19);
		\draw [style=nodirection] (19) to (18);
		\draw [style=nodirection] (18) to (17);
		\draw [style=nodirection] (16) to (19);
		\draw [style=nodirection] (16) to (20);
		\draw [style=nodirection] (17) to (23);
		\draw [style=nodirection] (19) to (22);
		\draw [style=nodirection] (18) to (21);
		\draw [style=nodirection] (30) to (26);
		\draw [style=nodirection] (26) to (27);
		\draw [style=nodirection] (27) to (24);
		\draw [style=nodirection] (29) to (31);
		\draw [style=nodirection] (31) to (25);
		\draw [style=nodirection] (25) to (28);
		\draw [style=nodirection] (28) to (29);
		\draw [style=nodirection] (29) to (24);
		\draw [style=nodirection] (24) to (30);
		\draw [style=nodirection] (24) to (25);
		\draw [style=nodirection] (25) to (27);
		\draw [style=nodirection] (27) to (31);
		\draw [style=nodirection] (31) to (26);
		\draw [style=nodirection] (26) to (29);
		\draw [style=nodirection] (29) to (30);
		\draw [style=nodirection] (30) to (28);
		\draw [style=nodirection] (28) to (24);
		\draw [style=nodirection] (24) to (26);
		\draw [style=nodirection] (27) to (30);
		\draw [style=nodirection] (29) to (25);
		\draw [style=nodirection] (28) to (31);
	\end{pgfonlayer}
\end{tikzpicture}\\
\hline
\end{tabular}
\caption{Graphs on $n$ vertices attaining the minimum trace with respect to $\amxL(X)$ for $n = 3,4,5,6,7,8$.\label{tab:Lmintr}}
\end{table}

\begin{table}[htp]
\centering
\begin{tabular}{|l|l|c|}
\hline
$n$ & $\min_X \tr(\amxA{X})$ & Graphs \\
\hline
$3$ & $1\nicefrac{1}{4}$ & \begin{tikzpicture}[scale=0.6]
	\begin{pgfonlayer}{nodelayer}
		\node [style=new] (0) at (0, 0) {};
		\node [style=new] (1) at (-1, 0) {};
		\node [style=new] (2) at (1, 0) {};
	\end{pgfonlayer}
	\begin{pgfonlayer}{edgelayer}
		\draw [style=nodirection] (1) to (0);
		\draw [style=nodirection] (0) to (2);
	\end{pgfonlayer}
\end{tikzpicture} \\
\hline
$4$ & $1\nicefrac{1}{5}$ & \begin{tikzpicture}[scale=0.6]
	\begin{pgfonlayer}{nodelayer}
		\node [style=new] (0) at (0, 0) {};
		\node [style=new] (1) at (-1, 0) {};
		\node [style=new] (2) at (1, 0) {};
		\node [style=new] (3) at (-2, 0) {};
	\end{pgfonlayer}
	\begin{pgfonlayer}{edgelayer}
		\draw [style=nodirection] (1) to (0);
		\draw [style=nodirection] (0) to (2);
		\draw [style=nodirection] (3) to (1);
	\end{pgfonlayer}
\end{tikzpicture} \\
\hline
$5$ & $1\nicefrac{1}{3}$ & \begin{tikzpicture}[scale=0.6]
	\begin{pgfonlayer}{nodelayer}
		\node [style=new] (0) at (0, 0) {};
		\node [style=new] (1) at (-1, 0) {};
		\node [style=new] (2) at (1, 0) {};
		\node [style=new] (3) at (-2, 0) {};
		\node [style=new] (4) at (-3, 0) {};

	\end{pgfonlayer}
	\begin{pgfonlayer}{edgelayer}
		\draw [style=nodirection] (1) to (0);
		\draw [style=nodirection] (0) to (2);
		\draw [style=nodirection] (3) to (1);
		\draw [style=nodirection] (4) to (3);
	\end{pgfonlayer}
\end{tikzpicture} \\
\hline
$6$ & $1\nicefrac{1}{4}$ & \begin{tikzpicture}[scale=0.6]
	\begin{pgfonlayer}{nodelayer}
		\node [style=new] (0) at (0, 1) {};
		\node [style=new] (1) at (0, 0) {};
		\node [style=new] (2) at (-1, 1) {};
		\node [style=new] (3) at (-1, 0) {};
		\node [style=new] (4) at (1, 1) {};
		\node [style=new] (5) at (1, 0) {};
	\end{pgfonlayer}
	\begin{pgfonlayer}{edgelayer}
		\draw [style=nodirection] (2) to (0);
		\draw [style=nodirection] (0) to (4);
		\draw [style=nodirection] (4) to (5);
		\draw [style=nodirection] (5) to (1);
		\draw [style=nodirection] (1) to (0);
		\draw [style=nodirection] (2) to (3);
		\draw [style=nodirection] (3) to (1);
	\end{pgfonlayer}
\end{tikzpicture}\\
\hline
$7$ & $1.349025083599055$ & \begin{tikzpicture}[scale=0.6]
	\begin{pgfonlayer}{nodelayer}
		\node [style=new] (0) at (0, 1) {};
		\node [style=new] (1) at (-1, 0.5) {};
		\node [style=new] (2) at (1, 0.5) {};
		\node [style=new] (3) at (-1.25, -0.5) {};
		\node [style=new] (4) at (1.25, -0.5) {};
		\node [style=new] (5) at (-0.5, -1.25) {};
		\node [style=new] (6) at (0.5, -1.25) {};
	\end{pgfonlayer}
	\begin{pgfonlayer}{edgelayer}
		\draw [style=nodirection] (0) to (2);
		\draw [style=nodirection] (2) to (4);
		\draw [style=nodirection] (4) to (6);
		\draw [style=nodirection] (5) to (3);
		\draw [style=nodirection] (3) to (1);
		\draw [style=nodirection] (1) to (0);
		\draw [style=nodirection] (0) to (4);
		\draw [style=nodirection] (4) to (5);
		\draw [style=nodirection] (5) to (1);
		\draw [style=nodirection] (1) to (2);
		\draw [style=nodirection] (2) to (6);
		\draw [style=nodirection] (6) to (3);
		\draw [style=nodirection] (3) to (0);
		\draw [style=nodirection] (3) to (4);
	\end{pgfonlayer}
\end{tikzpicture}\\
\hline
$8$ & $1\nicefrac{29}{185}$ & \begin{tikzpicture}[scale=0.6]
	\begin{pgfonlayer}{nodelayer}
		\node [style=new] (0) at (-3, 1) {};
		\node [style=new] (1) at (-3, -0.5) {};
		\node [style=new] (2) at (-2, 1.5) {};
		\node [style=new] (3) at (-2, 0) {};
		\node [style=new] (4) at (-1, 1) {};
		\node [style=new] (5) at (-1, -0.5) {};
		\node [style=new] (6) at (0, 1.5) {};
		\node [style=new] (7) at (0, 0) {};
	\end{pgfonlayer}
	\begin{pgfonlayer}{edgelayer}
		\draw [style=nodirection] (0) to (2);
		\draw [style=nodirection] (2) to (3);
		\draw [style=nodirection] (3) to (1);
		\draw [style=nodirection] (1) to (2);
		\draw [style=nodirection] (0) to (3);
		\draw [style=nodirection] (0) to (1);
		\draw [style=nodirection] (4) to (6);
		\draw [style=nodirection] (6) to (7);
		\draw [style=nodirection] (7) to (4);
		\draw [style=nodirection] (4) to (5);
		\draw [style=nodirection] (6) to (5);
		\draw [style=nodirection] (5) to (7);
		\draw [style=nodirection] (2) to (4);
		\draw [style=nodirection] (4) to (0);
		\draw [style=nodirection] (2) to (6);
		\draw [style=nodirection] (3) to (5);
		\draw [style=nodirection] (5) to (1);
		\draw [style=nodirection] (3) to (7);
	\end{pgfonlayer}
\end{tikzpicture} \\
\hline
\end{tabular}
\caption{Graphs on $n$ vertices attaining the minimum trace with respect to $\amxA$ for $n = 3,4,5,6,7,8$.\label{tab:Amintr}}
\end{table}

\section{Open problems}

In Section \ref{sec:constdiag}, we look graphs whose average mixing matrix has constant diagonal and give two constructions for such graphs; walk-regular graphs and their rooted products with $K_2$. It seems surprising that this property can occur for graphs which are not regular, but the computational results suggest that there would be other constructions for such graphs. Such constructions would be interesting. 

Based on the computations summarized in Table \ref{tab:Amaxtr} and on Corollary \ref{cor:regKnMaxTr}, we also make the following conjecture. 
\begin{conjecture}
The complete graph on $n$ vertices attains the maximum trace with respect to the $\mxm_{A}$ for all $n$. 
\end{conjecture}


\end{document}